\documentclass[letterpaper, 10 pt, conference]{ieeeconf}  

\IEEEoverridecommandlockouts                              

\usepackage{amsmath}
\usepackage{amssymb}

\usepackage{graphics}

\usepackage{color}

\usepackage{subcaption}

\usepackage{tikz}
\usetikzlibrary{arrows}
\usepackage{verbatim}

\usepackage{balance}

\newcommand{\dist}{{\rm dist}}
\newcommand{\cl}{{\rm cl}}

\newcommand{\mR}{{\mathbb R}}
\newcommand{\mC}{{\mathbb C}}

\newcommand{\Linf}{L_\infty}
\newcommand{\Hinf}{\mathcal{H}_\infty}
\newcommand{\HinfC}{\Hinf(\mC_+)}
\newcommand{\qed}{\hspace*{\fill}~$\blacksquare$\par}

\newtheorem{thm}{Theorem}

\newtheorem{lemma}[thm]{Lemma}

\newtheorem{remark}{Remark}

\newtheorem{ex}[thm]{Example}

\graphicspath{{Figures/}}

\title{\LARGE \bf
On analytic interpolation with non-classical constraints \\ for solving problems in robust control
}

\author{
Axel Ringh$^{1}$, Johan Karlsson$^{2}$, and Anders Lindquist$^{3,2}$
\thanks{*This work was supported by Knut and Alice Wallenberg foundation, grant KAW 2018.0349,  the Swedish Research Council (VR), grant 2014-5870, and SJTU-KTH cooperation grant.}
\thanks{$^{1}$Department of Electronic and Computer Engineering, The Hong Kong University of Science and Technology, Clear Water Bay, Kowloon, Hong Kong, China. {\tt\small eeringh@ust.hk}}
\thanks{$^{2}$Division of Optimization and Systems Theory, Department of Mathematics, KTH Royal Institute of Technology, Stockholm, Sweden. {\tt\small johan.karlsson@math.kth.se}}%
\thanks{$^{3}$Department of Automation and School of Mathematics, Shanghai Jiao Tong University, Shanghai, China. {\tt\small alq@math.kth.se}}
}

\begin{document}


\maketitle
\thispagestyle{empty}
\pagestyle{empty}

\setlength{\voffset}{4pt}


\begin{abstract}
In this work we consider robust stabilization of uncertain dynamical systems and show that this can be achieved by solving a non-classically constrained analytic interpolation problem.
In particular, this non-classical constraint confines the range of the interpolant, when evaluated on the imaginary axis, to a frequency-dependent set.
By considering a sufficient condition for when this interpolation problem has a solution, we derive an approximate solution algorithm that can also be used for controller synthesis.
Finally, the theory is illustrated on a numerical example with a plant with uncertain gain, phase, and output delay.
\end{abstract}

\section{Introduction}\label{sec:introduction}

The problem of robustly stabilizing an uncertain dynamical system is an important problem in systems and control.
Model uncertainties may stem from, e.g., unmodeled dynamics, the use of models obtained via system identification,
or the use of a reduced-order model for control design.
In fact, any mathematical model of a physical system will inevitably be incorrect to a certain degree.
Consequently, robust control has been the subject of much study; see, e.g., one of the monographs \cite{zhou1996robust, foias1996robust}, and also \cite[Chp.~4]{doyle1992feedback}, \cite[Chp.~12]{astrom2008feedback}.

For certain classes of both structured and unstructured disturbances, the robust stabilization problem can be understood in terms of the Nyquist plot and the Nyquist stability criterion, see, e.g., \cite{doyle1981multivariable, chen1982necessary}. More precisely, for a set of structured disturbances of the transfer function, a controller that stabilizes a nominal plant will robustly stabilize the plant if the Nyquist curve of the nominal feedback interconnection does not intersect certain regions in the complex plane; see, e.g. \cite[Sec.~9.5]{zhou1996robust}, \cite[Sec.~9.3]{astrom2008feedback}, \cite[pp.~51-52]{doyle1992feedback} for an example of this in terms of the gain margin problem. Similarly, for the unstructured case, it can be understood as stabilizing a ``ball of plants'' around the nominal Nyquist curve, cf. \cite[Sec.~9.1-9.3]{zhou1996robust},  \cite[Sec.~12.2]{astrom2008feedback} \cite[Sec.~4.3]{doyle1992feedback}.

Another way to formulate and solve robust stabilization problems has been to cast them as non-classically constrained analytic interpolation problems in terms of the sensitivity or complementary sensitivity function, see, e.g., \cite{tannenbaum1980feedback, zames1983feedback, kimura1984robust, khargonekar1985non, khargonekar1986robust, qi2017fundamental, ringh2018lower, ringh2019ananalytic}. The non-classical constraint in question is a restriction of the range of the interpolant to some set,  which could potentially vary with frequency. The problems that have  previously been solved are in general either sets that are
constant, i.e., frequency-independent \cite{tannenbaum1980feedback, khargonekar1985non}, or sets that correspond to discs with a radius that is frequency-dependent \cite{zames1983feedback, kimura1984robust, khargonekar1985non, khargonekar1986robust, qi2017fundamental}.

A particular robust stabilization problem that has received considerable attention lately is the maximum delay margin problem, see, e.g., \cite{middleton2007achievable, qi2014fundamental, qi2017fundamental, ju2016further, ringh2018lower, ringh2019ananalytic}. This is a robust stabilization problem with a structured, infinite-dimensional uncertainty.
The problem is easily understood in the Nyquist plot, and in \cite{ringh2018lower, ringh2019ananalytic} this was used to reformulate it as an analytic interpolation problem with a non-classical constraint, akin to the problems in \cite{tannenbaum1980feedback, zames1983feedback, kimura1984robust, khargonekar1985non, khargonekar1986robust}. However, in contrast to these results the structured disturbance in \cite{ringh2018lower, ringh2019ananalytic} gives rise to a frequency-dependent constraint that does not correspond to discs with varying radius. Moreover, the methods for unstructured disturbances in  \cite{zames1983feedback, kimura1984robust, khargonekar1985non} in general lead to estimates of the maximum delay margin that are conservative
\cite[Sec.~III and VIII.A]{ringh2019ananalytic}.
Therefore, in  \cite{ringh2018lower, ringh2019ananalytic}, new, less conservative estimation methods for the maximum delay margin were developed based on the non-classically constrained analytical interpolation formulation of the problem.

The above observations give rise to a number of questions.  When can constraints for robust stabilization that are formulated in terms of the Nyquist plot be mapped to non-classical constraints in the analytic interpolation problem? Which robust stabilization problems can be formulated as non-classically constrained analytic interpolation problems,  for which the restriction on the range of the interpolant may vary with frequency? Among these problems, which ones can be solved analytically?
For problems that cannot be solved analytically,
can one find good approximation methods?
In this work, we take steps towards an answer to some of these questions. 
More precisely, we extend the methods from \cite{ringh2019ananalytic}, developed for handling an uncertain time delay, to handle a larger class of structured, stable disturbances. In particular, we cast the corresponding robust stabilization problem as an analytic interpolation problem with a possibly frequency-dependent non-classical constraint on the range of the interpolant,
and thus
partly extends results in \cite{tannenbaum1980feedback, khargonekar1985non}.
Moreover, we also show how one could design approximate solution algorithms for such problems, similar to what was done for time delay systems in \cite{ringh2019ananalytic}.

The outline of the paper is as follows: in Section~\ref{sec:background} we introduce some background material and make the problem formulation precise. In Section~\ref{sec:interpolation_results} we show that robust stabilization can be achieved by solving a non-classically constrained analytic interpolation problem, and in Section~\ref{sec:approx_method} we introduce a sufficient condition for solvability and show how this can be used to derive approximate solution algorithms. Section~\ref{sec:num_ex} contains a numerical example where we derive a controller that robustly stabilize a plant with uncertain gain, phase, and output delay. Finally, Section~\ref{sec:conclusions} contain some concluding remarks. Some proofs are deferred to the appendix for improved readability.

\section{Background and problem formulation}\label{sec:background}
In this section we will present some background material and make the problem formulation precise.
Moreover, the section is also used to set up the notation used in the paper.

\subsection{Notation}
Let $\mC_+$ denote the open right half plane, let $\bar{\mC}_+ := \cl(\mC_+)$ denote the closure, and let $\tilde{\mC}_+$ denote the extended closure, i.e.,
$\tilde{\mC}_+ := \bar{\mC}_+ \cup \{ \infty\}$ (cf. \cite{youla1974single, tannenbaum1980feedback, khargonekar1985non, khargonekar1986robust}).
Moreover, let $\HinfC$ denote the space of functions that are bounded and analytic in $\mC_+$. This is a Banach space when equipped with the norm $\| F(s) \|_{\Hinf} := \sup_{s \in \mC_+} |F(s)| = \sup_{\omega \in \mR} | F(i\omega) | =: \| F(i\omega) \|_{\Linf}$, see, e.g., \cite{foias1996robust, hoffman1962banach}.
Furthermore, let $\dist(A,x) := \inf_{y \in A} | x - y |$ be the distance between a set $A$ and a point $x$; in particular, we define $\dist(\emptyset,x) = \infty$ for all $x \in \mC$, where $\emptyset$ denotes the empty set. Finally, for any element $x$ we define $x \not \in \emptyset$ to always be true.

\subsection{Problem formulation}
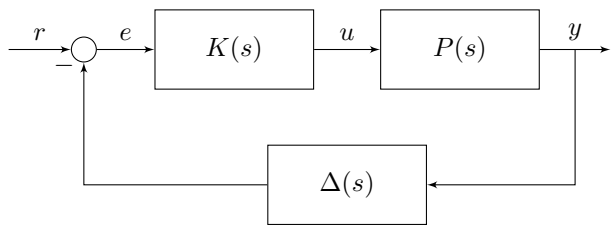
\begin{figure}[bt]
\tikzstyle{int}=[draw, minimum size=2em]
\tikzstyle{init} = [pin edge={to-,thin,black}]
\tikzstyle{block} = [draw, rectangle, 
    minimum height=3em, minimum width=6em]
\tikzstyle{sum} = [draw, circle, node distance=1cm]
\tikzstyle{input} = [coordinate]
\tikzstyle{output} = [coordinate]
\tikzstyle{pinstyle} = [pin edge={to-,thin,black}]

\begin{tikzpicture}[auto, node distance=2cm,>=latex']
    \node [input, name=input] {};
    \node [sum, right of=input] (sum) {};
    \node [block, right of=sum] (controller) {$K(s)$};
    \node [block, right of=controller, node distance=3cm] (system) {$P(s)$};

    \draw [->] (controller) -- node[name=u] {$u$} (system);
    \node [output, right of=system] (output) {};
    \node [block, below of=u] (delay) {$\Delta(s)$};

    \draw [draw,->] (input) -- node {$r$} (sum);
    \draw [->] (sum) -- node {$e$} (controller);
    \draw [->] (system) -- node [name=y] {$y$}(output);
    \draw [->] (y) |- (delay);
    \draw [->] (delay) -| node[pos=0.99] {$-$} (sum);
\end{tikzpicture}
\caption{
A feedback interconnection between a controller $K$, a plant $P$, and an uncertainty $\Delta$.}
\label{fig:blockdiagram}
\end{figure}

Let $P(s)$ be the transfer function of a
finite-dimensional, single-input-single-output (SISO) linear time-invariant (LTI) system, let $K(s)$ be a SISO LTI feedback controller,
and consider the feedback interconnection depicted in Fig.~\ref{fig:blockdiagram}. Here, $\Delta$ is a potential perturbation of the system. In this work, we are interested in the following robust stabilization problem: given a set of disturbances $\Omega$,
does there exist a controller $K$ that stabilizes the feedback interconnection in Fig.~\ref{fig:blockdiagram} for all $\Delta \in \Omega$?
If such a controller exists, we say that \emph{$K$ robustly stabilizes $P$ with respect to $\Omega$}.
The remaining problem is then how to find such a controller.

\subsection{Classical results for the unperturbed problem}
Before continuing we  review some classic results for the unperturbed stabilization problem, i.e., 
where $\Omega$ is the singleton
$\{ I \}$  with $I(s) \equiv 1$ being the one-function.
To this end, let $p_1, \ldots, p_n$ be the unstable poles and $z_1, \ldots, z_m$ the nonminimum phase zeros of $P$, respectively; for simplicity we will assume that they are all in $\mC_+$ and that they are all distinct.
In this case, $K$ is an (internally) stabilizing controller if and only if there is no pole-zero cancellation between $P$ and $K$ in $\tilde{\mC}_+$ and 
\begin{equation}
\label{eq:stability}
1 + P(s)K(s) \neq 0 \quad \text{for all } s \in \tilde{\mC}_+,
\end{equation}
see, e.g., \cite[p.~37]{doyle1992feedback}, \cite[p.~13]{helton1998classical}.
In particular, that it holds for $s = \infty$ means that there is no sequence $(s_n)_n \subset \bar{\mC}_+$ such that $\lim_{n\to \infty} 1 + P(s_n)K(s_n) = 0$.
These conditions can be equivalently expressed as an analytic interpolation problem in terms of the complementary sensitivity function $T := PK / (1 + PK)$, namely $T \in \HinfC$ and
\begin{subequations}\label{eq:interpolation}
\begin{align}
  &  T(p_k) = 1,\quad k = 1,\ldots, n , \\
  &  T(z_\ell)= 0,\quad \ell = 1,\ldots, m;
\end{align}
\end{subequations}
see, e.g., \cite{youla1974single}, \cite[Ch.~2 and 7]{helton1998classical}, \cite[pp.~72-73]{foias1996robust}.

For a $T \in \HinfC$ that satisfies \eqref{eq:interpolation}, a stabilizing controller for $P$ is given by $K = T / (P(1 - T))$.
This expression is easily obtained by algebraic manipulations. However, let us point out the fact 
that the open-loop transfer function $PK$ can be mapped to the complementary sensitivity function $T$ via the M\"obius transformation
\begin{equation}\label{eq:mobius}
\rho(s) = s/(1 + s).
\end{equation}
Moreover, the inverse mapping, which maps $T$ to $PK$ and which can be used to derive the form of the controller given above, is given by $\rho^{-1}(s) = s/(1 - s)$.

\begin{remark}\label{rem:inf_dim}
The assumption that $P$ is finite-dimensional can be relaxed to allow for certain infinite-dimensional systems, see \cite{callier1978algebra, desoer1980generalized, chen1982necessary, khargonekar1986robust} and \cite[pp.~56-57 and 72-73]{foias1996robust}. However, for the sake of simplicity we will derive and state the results in the finite-dimensional setting.
\end{remark}

\subsection{Assumptions on the set of perturbations $\Omega$}\label{sec:ass}
In order to give conditions for when robust stabilization is possible, we need to impose some restrictions on the set of perturbations $\Omega$. To start with, we
make the simplifying assumption that $\Omega \subset \HinfC \cap \mathcal{C}(i\mR)$, i.e., we consider stable disturbances that are continuous on the imaginary axis.
Moreover, we assume that $I \in \Omega$,
meaning that we consider stable disturbances around the nominal plant $P$. We also assume that no disturbance $\Delta$ cancels unstable poles in $P$, i.e., for all $\Delta \in \Omega$, 
$\Delta(p_k) \neq 0$ for $k = 1, \ldots, n$.%
\footnote{If there is a $\Delta \in \Omega$ for which $P\Delta$ have a pole-zero cancellation, we cannot find a controller $K$ that makes the system internally stable for this disturbance, and thus we cannot robustly stabilize $P$ with respect to $\Omega$. Moreover, note that no $\Delta$ can cancel nonminimum phase zeros of $P$, since $\Delta \in \HinfC$ and thus $|\Delta(z_\ell)| < \infty$ for $\ell = 1, \ldots, m$.}

We  also make some more technical assumptions: we assume that $\Omega$ is bounded,
that $\Omega$ is path-connected with respect to pointwise convergence, and
that for every bounded sequence $(\Delta_n)_n$ there exists a subsequence that converges pointwise in $\bar{\mC}_+$ to some $\Delta \in \Omega$.
In particular, bounded means that there is a common upper bound for the norm of all $\Delta \in \Omega$, i.e., $\sup_{\Delta \in \Omega} \| \Delta \|_{\Hinf} =: N < \infty$.
Moreover, 
path-connected with respect to pointwise convergence means that for any two elements $\Delta_1$, $\Delta_2$, we can find a curve $\gamma : [0,1] \to \Omega$ so that $\gamma(0) = \Delta_1$ and $\gamma(1) = \Delta_2$, and so that for all points $s \in \bar{\mC}_+$ the function $[\gamma(t)](s)$ is continuous in $t$.
Finally, by Montel's theorem (see, e.g., \cite[p.~35]{schiff1993normal}, \cite[Thm.~14.6]{rudin1987real})
the boundedness of $\Omega$ implies that $\Omega$ is a so called normal family, i.e., that any sequence of functions in $\Omega$ contains a subsequence that converges uniformly to an analytic function, call it $f$, on all compact subsets $A \subset \mC_+$.  
The last assumption, on pointwise convergence in $\bar{\mC}_+$ to an element in $\Omega$,  is an amplification of this to that the subsequence converges pointwise also on $i\mR$, and that $f \in \Omega$.

The reason for imposing these constraints is that when we vary $\Delta \in \Omega$ we want the corresponding Nyquist curve, i.e., the set $\{ \Delta(i\omega) \mid \omega \in \mR \}$, to remain bounded and to vary continuously pointwise. Moreover, since $I \in \Omega$ and $\Omega$ is path-connected, any disturbed plant $P\Delta$ can be obtained by pointwise continuously deforming the Nyquist curve of the nominal plant. 

\begin{ex}\label{ex:uncertain_zero}
As a motivating example, consider the SISO LTI system $(s + z)P$, where $P$ is a strictly proper plant. Assume that the exact position of the real minimum phase zero $z$ is uncertain, with upper and lower bounds $b > 0$ and $a > 0$, respectively.
To model this uncertainty as disturbances,
let $\Delta_y(s) = (s + y)/(s + z)$ for some $y \in [a,b]$ and consider the family of disturbances $\Omega = \{ \Delta_y(s) \mid y \in [a,b] \}$.
For this set of disturbances, most of our assumptions are easily verified, and thus we only consider the last two.
To this end,
to see that $\Omega$ is path connected with respect to pointwise convergence in $\bar{\mC}_+$, let $\Delta_{y_1}$ and $\Delta_{y_2}$ be two disturbances in $\Omega$.
Moreover, let $y(t)$ for $t \in [0,1]$ be a path connecting the real numbers $y_1$ and $y_2$, and in particular $y(t)$ is thus continuous in $t$. Now, let $\gamma(t) = (s + y(t))/(s + z)$ and note that $[\gamma(t)](s)$ is continuous in $t$ for all $s \in \bar{\mC}_+$. Therefore, $\gamma$ is a path connecting $\Delta_{y_1}$ and $\Delta_{y_2}$, and hence the set is path connected. 
Next, let $(\Delta_{y_n}(s))_n \subset \Omega$ be any sequence of disturbances. Since $(y_n)_n \subset [a, b]$ is a bounded sequence of real numbers, there is a converging subsequence $(y_{n_k})_{n_k}$ such that $y_{n_k} \to y_{\infty} \in [a, b]$. In particular, this means that $\Delta_{y_{n_k}}(s) \to \Delta_{y_{\infty}}(s) \in \Omega$ pointwise in $\bar{\mC}_+$.
\end{ex}

Note that a system with an uncertainty in a stable pole can be handled analogously to the above example. In fact, a system with simultaneous uncertainties in minimum phase zeros, stable poles, gain, phase, and delay in the feedback loop, can be handled analogously.

\subsection{A Nyquist-type result for robust stabilization}
As a preliminary, we first derive a Nyquist-type result for robust stabilization, cf. \cite{doyle1981multivariable, chen1982necessary}. It is derived under the simplifying assumption that $P$ is strictly proper. However
this is only done in order to not clutter the main idea with technical details.
To this end, note that for $K$ to stabilize the feedback loop in Fig.~\ref{fig:blockdiagram} for all $\Delta \in \Omega$, we need to alter \eqref{eq:stability} to read
\begin{equation}
\label{eq:robust_stability}
1 + P(s)K(s)\Delta(s) \neq 0 \text{ for all } s \in \tilde{\mC}_+ \text{ and all } \Delta \in \Omega.
\end{equation}
Moreover, we need that no pole-zero cancellation occurs between $K$ and any of the plants $P\Delta$.

Next, consider the set-valued function
\begin{equation}\label{eq:nyquist_set_valued_function}
\Gamma(i\omega) := \left\{ z  \; \Big| \; z = \frac{-1}{\Delta(i\omega)}, \; \Delta \in \Omega \text{ and } \Delta(i\omega) \neq 0 \right\},
\end{equation}
where, for a given frequency $\omega$, we let $\Gamma(i\omega) = \emptyset$ if $\Delta(i\omega) = 0$ for all $ \Delta \in \Omega$.
Now, note that given an $\omega_0$, for each $z \in \Gamma(i\omega_0)$ there exists a $\Delta \in \Omega$ such that $z \Delta(i \omega_0) = -1$. Based on this observation, we have the following result.

\begin{lemma}\label{lem:robust_nyquist}
Assume that $P$ is strictly proper and that $\Omega$ satisfies the assumptions in Section~\ref{sec:ass}, and let $\Gamma$ be given as in \eqref{eq:nyquist_set_valued_function}.
Then a proper $K$ robustly stabilizes $P$ with respect to $\Omega$ if and only if $K$ stabilizes $P$ and $K(i\omega)P(i\omega) \not \in \Gamma(i\omega)$ for all $\omega \in \mR$.
\end{lemma}

\hspace{-22pt}
\begin{proof}
Without loss of generality, assume that $K$ has $n_K$ poles in $\mC_+$. Then
$K$ stabilizes $P$ if and only if the Nyquist curve $P(i\omega)K(i\omega)$ encircles $-1$ exactly $n + n_K$ times counter-clockwise without intersecting it, see, e.g., \cite[Thm.~9.2]{astrom2008feedback}, \cite[p.~37]{doyle1992feedback}. 
By the pointwise continuity of the Nyquist curve with respect to the disturbance, the result now follows since the number of encirclements and intersections can change if and only if at least one of the following statements is true: i) there exists a $\Delta$ and an $\omega_0$
such that
$K(i\omega_0)P(i\omega_0)\Delta(i\omega_0) = -1$, or ii) there exists a sequence $(\omega_n, \Delta_n)_n$, $|\omega_n| \to \infty$, such that $-1 = \lim_{n \to \infty} K(i\omega_n)P(i\omega_n)\Delta_n(i\omega_n)$.  ii) cannot happen since $K(i\omega_n)$ and $\Delta_n(i\omega_n)$ remain bounded while $P(i\omega_n)$ goes to 0, and by construction, i) cannot happen if and only if $K(i\omega)P(i\omega) \not \in \Gamma(i\omega)$ for all $\omega \in \mR$.
\end{proof}

\section{Non-classically constrained analytic interpolation for robust stabilization}\label{sec:interpolation_results}
The constraint on the Nyquist curve in Lemma~\ref{lem:robust_nyquist} can be reformulate as a non-classical constraint on the analytic interpolant $T$ via the M\"obius transformation \eqref{eq:mobius}.
To this end, since $I \in \Omega$, for $K$ to robustly stabilize $P$ with respect to $\Omega$ we have that $K$ must stabilize $P$. Therefore, we can restrict our attention to controllers $K$ such that \eqref{eq:stability} holds. 
For such controllers, adding and subtracting $PK$ on the left hand side in \eqref{eq:robust_stability} and dividing both sides by $1 + PK$ gives the condition
\begin{equation}\label{eq:nec_and_suff_cond}
1+ T(s)(\Delta(s) - 1) \neq 0 \text{ for all } s\in \tilde{\mC}_+ \text{ and all } \Delta \in \Omega.
\end{equation}
Again, that this holds for $s = \infty$ means that there is no sequence $(s_n)_n \subset \bar{\mC}_+$ such that $\lim_{n \to \infty} 1+ T(s_n)(\Delta(s_n) - 1) = 0$. Therefore,  \eqref{eq:nec_and_suff_cond} is in fact equivalent to the condition that for each $\Delta \in \Omega$ there is an $\varepsilon > 0$ such that
\begin{equation}\label{eq:nec_and_suff_cond_2}
|1+ T(s)(\Delta(s) - 1)| \geq  \varepsilon
\end{equation}
for all $s\in \bar{\mC}_+$.
Using this, we can prove the following intermediate result.
\begin{lemma}\label{lem:robust_interoplant}
Assume that $\Omega$ satisfies the assumptions in Section~\ref{sec:ass}. Then there exists a controller $K$ that robustly stabilizes $P$ with respect to $\Omega$ if and only if there exists a $T \in \HinfC$ that satisfies \eqref{eq:interpolation} and \eqref{eq:nec_and_suff_cond_2}.
\end{lemma}

\hspace{-22pt}
\begin{proof}
A controller $K$ robustly stabilizes $P$ with respect to $\Omega$ if and only if \eqref{eq:robust_stability} holds, and there is no pole-zero cancellation between $K$ and $P\Delta$ for all $\Delta \in \Omega$. By the arguments leading up to the lemma, this implies that there exists a $T \in \HinfC$ that satisfies \eqref{eq:interpolation} and \eqref{eq:nec_and_suff_cond_2}.

To prove the converse, assume such a $T$ exists and let $K$ be the corresponding stabilizing controller of $P$, i.e., let $T = PK/(1 + PK)$.
If $s_0 \in \bar{\mC}_+$ is a pole of $K$, then $T(s_0) = 1$. Assume that it is also a zero of $\Delta$; this gives $|1 + T(s_0)(\Delta(s_0) - 1)| = 0$, which contradicts \eqref{eq:nec_and_suff_cond_2}.   
Left to show is that $K$ and $\Delta$ cannot have a pole-zero cancellation in $s = \infty$. To this end, assume that they have a pole-zero cancellation in $s = \infty$ and let $(s_n)_n \subset \bar{\mC}_+$ be a sequence such that $|s_n| \to \infty$ as $n \to \infty$.
In particular, that $K$ has  a pole in $s = \infty$ means that $\lim_{n \to \infty} T(s_n) = 1$.
Therefore, 
\[
\lim_{n \to \infty} |1 + T(s_n)(\Delta(s_n) - 1)| = |1 + 1(0 -1)| = 0,
\] 
which contradicts \eqref{eq:nec_and_suff_cond_2}. Thus the result follows.
\end{proof}

Next, we show that \eqref{eq:nec_and_suff_cond_2} is equivalent to a non-classical constraint on $T$ on the imaginary axis, akin to the result in Lemma~\ref{lem:robust_nyquist}.
More specifically, consider
\begin{align}
\Lambda(i\omega) \!  := \!  \left\{ \! z \, \Big| \, z = \frac{1}{1 \! - \! \Delta(i\omega)}, \; \Delta \in \Omega \text{ and } \Delta(i\omega) \! \neq \! 1  \! \right\} \!, \label{eq:lambda}
\end{align}
which can be interpreted as $\Lambda(i\omega)  = \rho \circ \Gamma(i\omega)$.
To be precise, we also make the following definition:
\begin{align*}
\Lambda(\infty) := \cl \Big( \big\{ & z \in \mC \;  \big| \; z  = \lim_{n \to \infty} \frac{1}{1 - \Delta_n(s_n)}, \\
&  \text{ for }  (s_n, \Delta_n)_n \subset \bar{\mC}_+ \times \Omega, \; |s_n| \to \infty \big\} \Big). 
\end{align*}

In light of Lemma~\ref{lem:robust_nyquist} and the M\"obius transformation \eqref{eq:mobius},  a non-classical constraint on the analytic interpolant of the form $T(i\omega) \not \in \Lambda(i\omega)$ should be equivalent to \eqref{eq:nec_and_suff_cond_2}. That this is indeed true (in certain cases) is the result of the next theorem.

\begin{thm}\label{thm:dist_to_set}
Assume that $\Omega$ satisfies the assumptions in Section~\ref{sec:ass}, and let $\Lambda$ be given as in \eqref{eq:lambda}.
Let $T \in \HinfC \cap \mathcal{C}(i\mR)$, and assume that $T(\infty)$ is well-defined and that $T(\infty) \not \in \Lambda(\infty)$. Then the following two statements are equivalent:
\begin{enumerate}
\item there exists an $\varepsilon > 0$ such that $|1 + T(s)(\Delta(s) - 1)| \geq \varepsilon$ for all $s \in \tilde{\mC}_+$ and $\Delta \in \Omega$
\item there exists an $\varepsilon > 0$ such that $\dist ( \Lambda(i\omega), T(i\omega)) \geq \varepsilon$ for all $\omega \in \mR$.
\end{enumerate}
\end{thm}

\hspace{-22pt}
\begin{proof}
See Appendix~\ref{app:proof_of_thm:thm:dist_to_set}.
\end{proof}

From Lemma~\ref{lem:robust_interoplant} and Theorem~\ref{thm:dist_to_set}, we have the following result for robust stabilization of an uncertainty set $\Omega$.

\begin{thm}\label{cor:main_result}
Assume that $\Omega$ satisfies the assumptions in Section~\ref{sec:ass}. Then there exists a controller $K$ that robustly stabilizes $P$ with respect to $\Omega$ if
there is a solution to the non-classically constrained analytic interpolation problem:
find $T \in \HinfC \cap \mathcal{C}(i\mR)$ that satisfies
\begin{enumerate}
\item the interpolation conditions \eqref{eq:interpolation},
\item $T(i\omega) \not \in \Lambda(i\omega)$ for all $\omega \in \mR$,
\item $T(\infty) \not \in \Lambda(\infty)$.
\end{enumerate}
\end{thm}

\begin{ex}
We return to the example with an uncertain position of a minimum phase zero, as presented in Example~\ref{ex:uncertain_zero}. A direct calculation gives that the forbidden regions for the Nyquist plot and the for the analytic interpolant are
\begin{align*}
& \Gamma(i\omega) = 
\left\{- \frac{zy + \omega^2 + i \omega (y-z)}{y^2 + \omega^2} \; \Big| \; y \in [a,b]  \right\}, \\
& \Lambda(i\omega) = \left\{ \frac{z + i\omega}{z - y}  \; \Big| \; y \in [a,b] \right\}.
\end{align*}
\end{ex}

\section{Sufficient conditions for stabilization and a controller synthesis algorithm}\label{sec:approx_method}
Analytic interpolation problems with a non-classical constrain of the type given in Theorem~\ref{cor:main_result} have been solved in certain cases, see Section~\ref{sec:introduction}. However, to the best of our knowledge no general solvability conditions exist for such problems. In this section we will derive sufficient conditions for when such problems can be solved. Moreover, when these conditions are fulfilled, we also obtain an algorithm for constructing robustly stabilizing controllers.

To this end, we return to the conditions given in Lemma~\ref{lem:robust_interoplant} and note in particular that \eqref{eq:nec_and_suff_cond} is fulfilled if $\sup_{\Delta \in \Omega} \|T(s) (\Delta(s) - 1)  \|_{\HinfC} < 1$.  In fact, this can be interpreted as a small-gain argument on the systems $T$ and $\Delta - 1$, and from this we get that a sufficient condition for robust stabilization is that
\begin{equation}\label{eq:suff_cond}
\inf_{\substack{T \in \HinfC \\ \text{subject to } \eqref{eq:interpolation}}} \sup_{\Delta \in \Omega} \|T(s) (\Delta(s) - 1)  \|_{\HinfC} < 1.
\end{equation}
Since $T(s) (\Delta(s) - 1) \in \HinfC$, this can be rewritten as
\begin{align}\label{eq:suff_cond_reworked}
1 & >
\inf_{\substack{T \in \HinfC \\ \text{subject to } \eqref{eq:interpolation}}} \sup_{\Delta \in \Omega} \|T(i\omega) (\Delta(i\omega) - 1)  \|_{\Linf(i\mR)} \nonumber  \\
& = \inf_{\substack{T \in \HinfC \\ \text{subject to } \eqref{eq:interpolation}}} \sup_{\omega \in \mR} \left[ |T(i\omega)| \, \sup_{\Delta \in \Omega} |\Delta(i\omega) - 1| \right] \nonumber  \\
& = \inf_{\substack{T \in \HinfC \\ \text{subject to } \eqref{eq:interpolation}}} \sup_{\omega \in \mR} \left[ |T(i\omega)| \frac{1}{\inf_{\Delta \in \Omega} \frac{1}{|\Delta(i\omega) - 1|}} \right] \nonumber  \\
& = \inf_{\substack{T \in \HinfC \\ \text{subject to } \eqref{eq:interpolation}}} \sup_{\omega \in \mR} \big[ |T(i\omega)| \phi(i\omega) \big],
\end{align}
where 
\[
\phi(i\omega) :=  1 / \dist(\Lambda(i\omega), 0) = \sup_{\Delta \in \Omega} \! |\Delta(i\omega) \! - \! 1|.
\]
In \eqref{eq:suff_cond_reworked}, the first
equality follows by noticing that the norm is defined by a supremum and changing order on the two suprema, and the
second equality follows since $\inf_x 1/f(x) = 1/\sup_x f(x)$ for nonnegative functions $f$.

This shows that a small-gain approach to
\eqref{eq:nec_and_suff_cond}
confines the range of the interpolant $T$ to a disc centered in the origin, whose radius is inversely proportional to the distance to the closest point in $\Lambda(i\omega)$. Moreover, the weight function $\phi$ exhibits some desirable regularity properties.

\begin{lemma}\label{lem:regularity_phi}
Assume that $\Omega$ satisfies the assumptions in Section~\ref{sec:ass}. Then $\phi(i\omega)$ is a lower semicontinuous function, bounded between $0$ and $N +1$.
\end{lemma}

\hspace{-22pt}
\begin{proof}
See Appendix~\ref{app:proof_lem:regularity_phi}.
\end{proof}

The regularity properties of $\phi$ can be used to devise a numerical approximation algorithm
for investigating feasibility of \eqref{eq:suff_cond_reworked}.
To this end,
take $\epsilon > 0$ and consider $\phi^{\epsilon}(i\omega) := \max\{ \phi(i\omega), \epsilon \}$. Since $0 < \epsilon \leq \phi^{\epsilon}(i\omega) \leq N+1 < \infty$ for all $\omega$, we have that $\log(\phi^{\epsilon}(i\omega))/(1 + \omega^2)$ is integrable.%
\footnote{Note that a lower semicontinuous function is measurable.}
Therefore, $\phi^{\epsilon}$ can be interpreted as the magnitude of an outer function $W^{\epsilon} \in \HinfC$, i.e., $|W^{\epsilon}(i\omega)| =\phi^{\epsilon}(i\omega)$ almost everywhere on $i\mR$, where the values of $W^{\epsilon}(s)$ for $s \in \mC_+$ are given by \cite[pp.~132-133]{hoffman1962banach}
\begin{equation}
\label{eq:outerrepr}
W^{\epsilon}(s) = \exp \! \left[ \frac{1}{\pi} \int_{-\infty}^\infty \!\! \log \big( \phi^{\epsilon}(i\omega) \big) \frac{ \omega s \! + \! i}{\omega \! + \! is} \frac{1}{1 \! + \!\omega^2} \,d\omega \right] \! . \!
\end{equation}
Thus, if $ \phi(i\omega)$ can be computed efficiently, by numerical integration one can evaluate $W^{\epsilon}(s)$ in points $s \in \mC_+$.  Moreover, 
since $\phi^{\epsilon}(i\omega) \geq \phi(i\omega)$ for all $\omega$, we get that a sufficient condition for \eqref{eq:suff_cond_reworked} is that
\begin{equation*}
1 > \inf_{\substack{T \in \HinfC \\ \text{subject to } \eqref{eq:interpolation}}} \|T(s) W^{\epsilon}(s) \|_{\HinfC}
\end{equation*}
which is a weighted minimization of the complementary sensitivity function akin to \cite{zames1983feedback, kimura1984robust, khargonekar1985non, khargonekar1986robust},  cf. \cite[p.~4]{ringh2019ananalytic}.
Since $W^{\epsilon}$ is outer and strictly positive on $i\mR$,
feasibility of  \eqref{eq:suff_cond_reworked} can be investigated
by introducing the ``change of variable'' $\tilde{T} = TW^{\epsilon}$
and consider the equivalent feasibility problem: find $\tilde{T} \in \HinfC$ such that 
\begin{subequations} \label{eq:T_tilde}
\begin{align}
& \|\tilde{T}(s)\|_{\HinfC} < 1 \label{eq:T_tilde_a}  \\
& \tilde{T}(p_k) =  W^{\epsilon}(p_k), & k = 1, \ldots, n,  \label{eq:T_tilde_b}\\
& \tilde{T}(z_\ell) =  0, & \ell = 1, \ldots, m.  \label{eq:T_tilde_c} 
\end{align}
\end{subequations}
The latter is a standard Nevanlinna-Pick interpolation problem, where the values in the interpolation points can be evaluated by numerical integration. Moreover, \eqref{eq:T_tilde} has a solution if and only if the corresponding Pick matrix is positive definite, see, e.g., \cite[pp.~157-159]{doyle1992feedback}, \cite[Sec.~2.9]{foias1996robust}.

\subsection{Introducing a shift for a tunable method}\label{sec:approx_method_shift}
Based on the insights from the previous section we now introduce a shift in the interpolant, which gives a tunable method. To this end, let $T(s) = \hat{T}(s) + T_0(s)$, where $ \hat{T} \in \HinfC$ and $T_0 \in \HinfC \cap \mathcal{C}(i\mR)$. Then \eqref{eq:nec_and_suff_cond} can be rewritten as
\begin{align*}
\hat{T}(s)(\Delta(s) - 1) \neq -1 + T_0(s) (1 - \Delta(s)) &\\
\text{for all } s \in \tilde{\mC}_+ \text{ and all } \Delta \in \Omega.&
\end{align*}
By Theorem~\ref{thm:dist_to_set} we know that the right hand side is bounded away from zero if and only if $\dist(\Lambda(i\omega), T_0(i\omega))$ is bounded away from zero and $T_0(\infty) \not \in \Lambda(\infty)$. For such $T_0$, we can rewrite the condition to read:  for all $s \in \tilde{\mC}_+$ and all $\Delta \in \Omega$,
\[
\hat{T}(s) \frac{\Delta(s) - 1}{1 - T_0(s) (1 - \Delta(s))} \neq -1.
\]
Using the same type of small-gain argument as in \eqref{eq:suff_cond}, and repeating the arguments in \eqref{eq:suff_cond_reworked}, we get that a sufficient condition for robust stabilization is that 
\begin{align}\label{eq:suff_cond_shifted}
1 & > \inf_{\substack{\hat{T} \in \HinfC \\ \text{subject to } \eqref{eq:interpolation}}} \sup_{\Delta \in \Omega} \left\|\hat{T}(s)  \frac{\Delta(s) - 1}{1 - T_0(s) (1 - \Delta(s))}  \right\|_{\HinfC} \nonumber \\
& = \inf_{\substack{\hat{T} \in \HinfC \\ \text{subject to } \eqref{eq:interpolation}}} \sup_{\omega \in \mR} \left[ |\hat{T}(i\omega)| \frac{1}{\inf_{\Delta \in \Omega} \left| T_0(i\omega) - \frac{1}{1 - \Delta(i\omega)} \right|} \right] \nonumber  \\
& = \inf_{\substack{T \in \HinfC \\ \text{subject to } \eqref{eq:interpolation}}} \sup_{\omega \in \mR} \big[ |\hat{T}(i\omega)| \phi_{T_0}(i\omega) \big],
\end{align}
where 
\[
\phi_{T_0}(i\omega) := \frac{1}{\dist(\Lambda(i\omega), T_0(i\omega))}.
\]
Therefore, $T_0(i\omega)$ can be interpreted as a shift of the center of the disc to which $\hat{T}(i\omega)$ is confined by the small-gain approach.
Moreover, since $\dist(\Lambda(i\omega), T_0(i\omega))$ is bounded from below, $\phi_{T_0}(i\omega)$ is bounded from above, and since $|(\Delta(i\omega) - 1)/(1 - T_0(i\omega) (1 - \Delta(i\omega)))|$ is continuous on $i\mR$, corresponding regularity properties as in Lemma~\ref{lem:regularity_phi} also hold for $\phi_{T_0}$. Therefore, introducing an $\epsilon > 0$ and considering $\phi^{\epsilon}_{T_0}(i\omega) := \max\{ \phi_{T_0}(i\omega), \epsilon \}$, we can ``extend'' $\phi^{\epsilon}_{T_0}(i\omega)$ to an outer function $W^{\epsilon}_{T_0}(s)$ via \eqref{eq:outerrepr}. Finally, with the ``change of variable'' $\tilde{T} = \hat{T}W^{\epsilon}_{T_0}$, feasibility of \eqref{eq:suff_cond_shifted} is equivalent to solvability of a Nevanlinne-Pick problem similar to \eqref{eq:T_tilde},
namely: find $\tilde{T} \in \HinfC$ such that
\begin{subequations} \label{eq:T_tilde_shifted}
\begin{align}
& \|\tilde{T}(s)\|_{\HinfC} < 1  \label{eq:T_tilde_shifted_a}  \\
& \tilde{T}(p_k) =  (1 - T_0(p_k))W^{\epsilon}_{T_0}(p_k), & k = 1, \ldots, n,  \label{eq:T_tilde_shifted_b}\\
& \tilde{T}(z_\ell) =   -T_0(z_\ell)W^{\epsilon}_{T_0}(z_\ell), &  \ell = 1, \ldots, m.  \label{eq:T_tilde_shifted_c} 
\end{align}
\end{subequations}

\begin{remark}
Note that although we must select $T_0$ such that $\dist(\Lambda(i\omega), T_0(i\omega))$ is bounded away from zero, $T_0$ need not fulfill the interpolation conditions \eqref{eq:interpolation}. On the other hand, if $T_0$ also fulfills the interpolation conditions, then $\tilde{T}(s) \equiv 0$ trivially fulfills \eqref{eq:T_tilde_shifted}. Based on this we can conclude that if the non-classically constrained interpolation problem in Theorem~\ref{cor:main_result} has a solution, this solution can be obtained as a solution to \eqref{eq:T_tilde_shifted} by a suitable choice of $T_0$.
\end{remark}

\subsection{Rational approximation for implementable controller}\label{sec:rat_approx}
To conclude this section, let us make a short remark on the controller synthesis problem.
To this end, let $\tilde{T} \in \HinfC$ such that $\| \tilde{T} \|_{\HinfC} < 1$, let $T_0 \in \HinfC \cap \mathcal{C}(i\mR)$ such that
$\dist(\Lambda(i\omega), T_0(i\omega))$ is bounded away from zero and $T_0(\infty) \not \in \Lambda(\infty)$ , and let $\Tilde{T}, T_0$ be such that they fulfill \eqref{eq:T_tilde_shifted_b} and \eqref{eq:T_tilde_shifted_c}. Then, a controller that robustly stabilizes $P$ with respect to $\Omega$ is given by
\[
K = \frac{\tilde{T} + T_0W^{\epsilon}_{T_0}}{P((1 - T_0) W^{\epsilon}_{T_0} - \tilde{T})}.
\]
However, even if $\tilde{T}$ and $T_0$ would be rational, $W^{\epsilon}_{T_0}$ is typically not rational, and thus $K$ is typically not rational. To overcome this, we can make a rational (over) approximation of $W^{\epsilon}_{T_0}$ and use that in the procedure. Moreover, by making an improper approximation of the weight, we can also guarantee that the obtained controller is proper. For the sake of brevity, we refer the reader to \cite[Sec.~VI.B-D]{ringh2019ananalytic}.

\section{Numerical example: controller design under uncertain gain, phase, and time delay}\label{sec:num_ex}
In this section we will apply the theory developed above to a numerical example to derive a rational, proper controller which robustly stabilize a given plant in the presence of uncertainty in the gain and the phase of the plant, and an uncertain time delay in the output signal.  More specifically, in this case the set of disturbances is given by
\[
\Omega = \left\{ \kappa e^{-i\theta} e^{-s\tau} \mid \kappa \in [1, k], \; \theta \in [-\varphi, \varphi], \; \tau \in [0, \mathtt{t}] \right\}.
\]
For the Nyquist curve $PK$, this gives the forbidden region
\[
\Gamma(i\omega) := \left\{ -\kappa^{-1} e^{i\theta} e^{i\omega\tau}  \mid \kappa \! \in \! [1, k], \theta \! \in \! [-\varphi, \varphi], \tau \! \in \! [0, \mathtt{t}]  \right\} \!.
\]
In particular, for $|\omega| \geq (2\pi - 2\varphi)/\mathtt{t}$ this corresponds to confining the Nyquist curve to a disc of radius $k^{-1}$. Therefore, any encirclement of $-1$ must be done in a relatively small window of frequencies. Similarly, for
the complementary sensitivity function this results in the frequency-dependent forbidden region 
\begin{align*}
& \Lambda(i\omega) \! = \! \left\{\! \frac{1}{1 \! - \! \kappa e^{-i\theta} e^{-i\omega\tau}}    \, \Big| \,    \kappa \! \in \! [1, k], \theta \! \in \! [-\varphi, \varphi],  \tau \! \in \! [0, \mathtt{t}] \! \right\} \!,
\end{align*}
and for $|\omega| \geq (2\pi - 2\varphi)/\tau$ we have that $\mC \setminus \Lambda(i\omega)$ is the disc $\{ z \mid |z + 1/(k^2 - 1)| < k/(k^2 -1) \}$, cf. \cite[App.~A]{ringh2019ananalytic}.

\begin{figure*}[tb]
\begin{center}
\begin{subfigure}[t]{.325\textwidth}
  \centering
  \includegraphics[width=\textwidth]{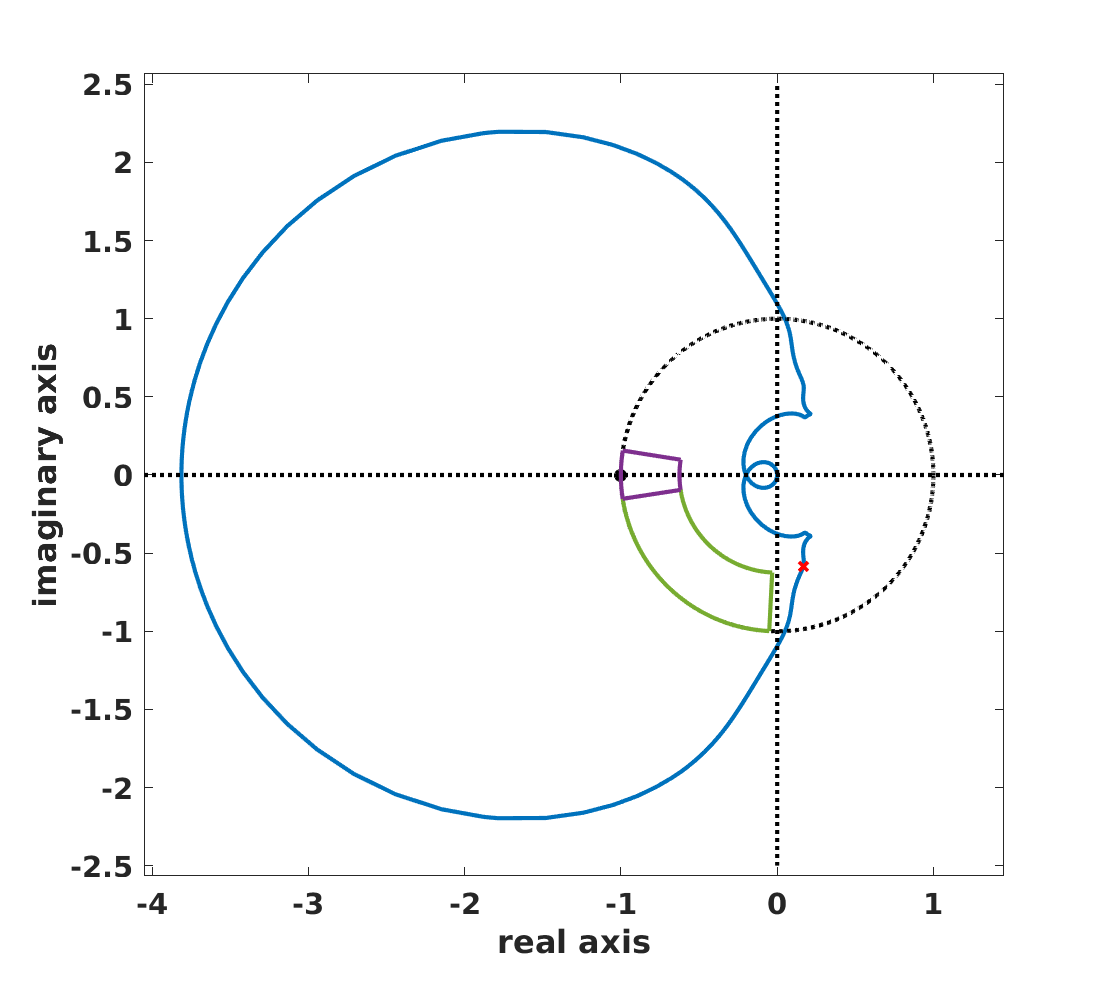}
  \subcaption{Frequency $\omega = 0.9138$.}
  \label{subfig:num_ex_nyquist_low_freq}
\end{subfigure}
\hfill
\begin{subfigure}[t]{.325\textwidth}
 \centering
  \includegraphics[width=\textwidth]{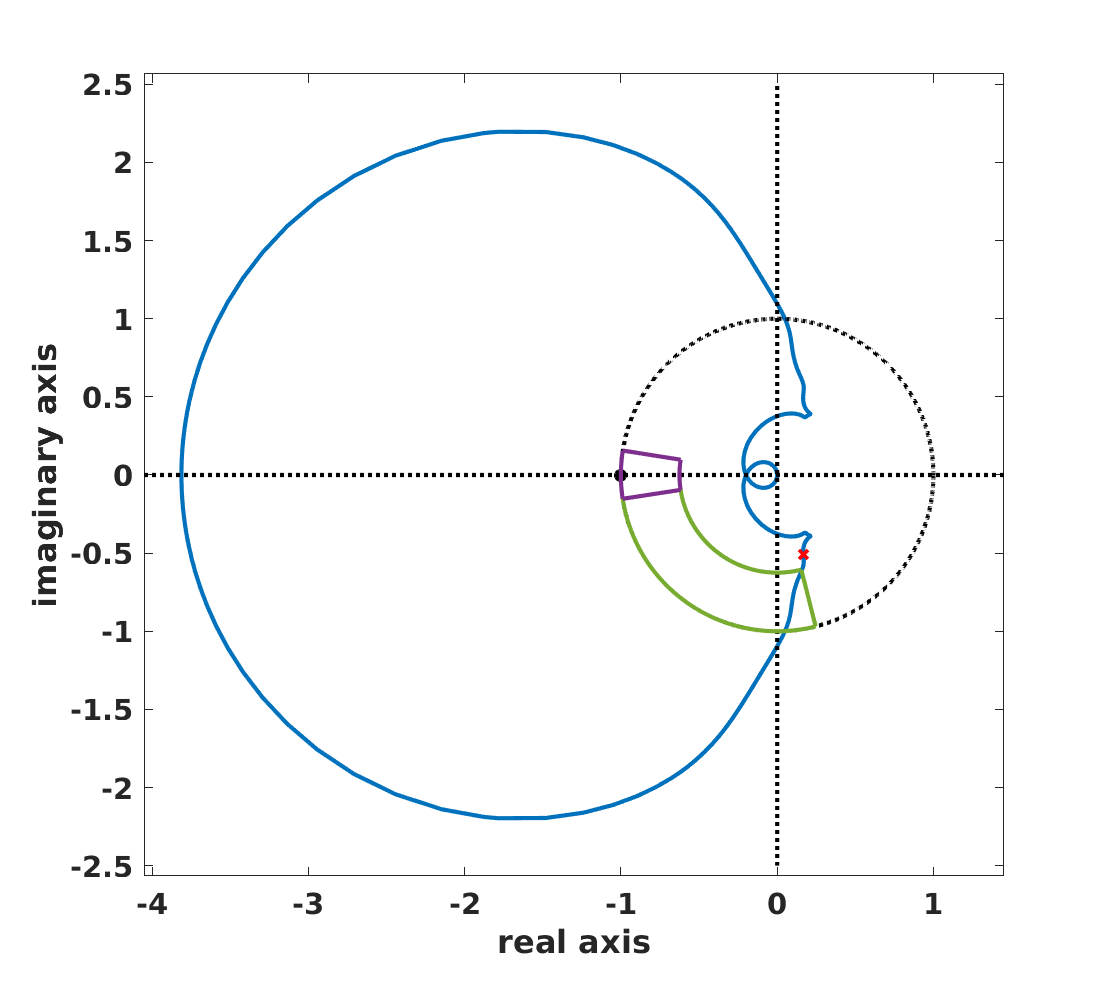}
  \subcaption{Frequency $\omega = 1.110$.}
  \label{subfig:num_ex_nyquist_med_freq}
\end{subfigure}
\hfill
\begin{subfigure}[t]{.325\textwidth}
  \centering
  \includegraphics[width=\columnwidth]{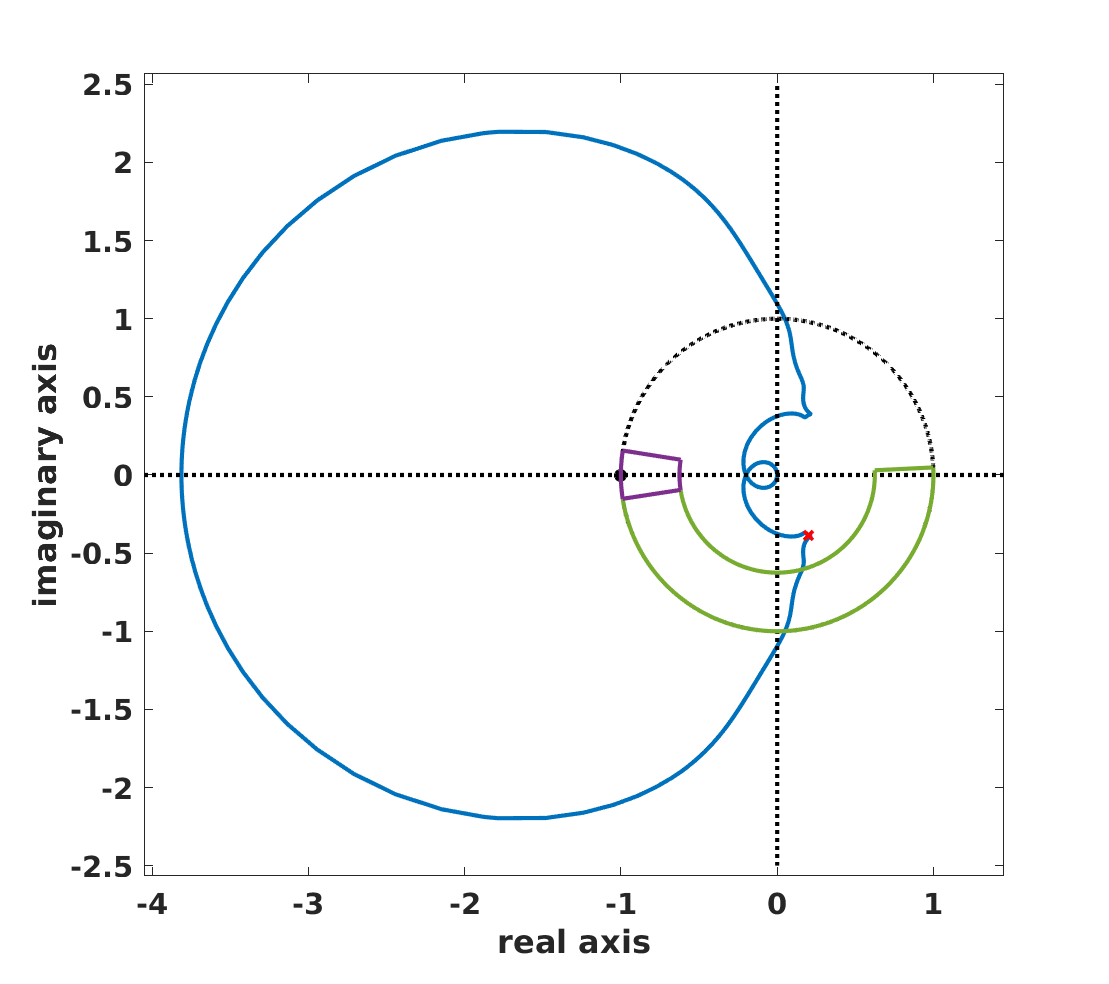}
  \caption{Frequency $\omega = 2.024$.}
  \label{subfig:num_ex_nyquist_high_freq}
\end{subfigure}
\caption{Nyquist plots of the obtained open-loop system \eqref{eq:PK_controller_design_ex}. The region outlined by the purple line segments corresponds to the constant set $\Gamma$ for a simultaneous gain and phase uncertainty in the intervals $[1, 1.6]$ and $[-\pi/20, \pi/20]$, respectively. Furthermore, the region outlined by the green line segments corresponds to the extension of the forbidden region obtained by adding the simultaneous delay uncertainty in the interval $[0, 1.5]$. Note that this set is frequency dependent, and in the three subfigures $\Gamma(i\omega)$ is illustrated for three different frequencies $\omega < (2\pi - 2\varphi)/\mathtt{t} = 19\pi/15 \approx 3.979$. For reference,  in each subfigure the corresponding point on the Nyquist plot is marked in red.}
\label{fig:num_ex_nyquist}
\end{center}
\end{figure*}

The example taken here is a modification of the example in \cite[Sec.~VIII.E]{ringh2019ananalytic}.%
\footnote{In particular, we use the same parameter settings in the approximation algorithm as in  \cite[Sec.~VIII]{ringh2019ananalytic}, cf. Section~\ref{sec:rat_approx}. Due to space limitations, we do not restate them here.}
To this end, consider the plant
\begin{equation}\label{eq:ex_plant}
P(s) =0.1 \frac{(0.1s - 1)(s + 0.1659)}{(s - 0.1081)(s^2 + 0.2981s + 0.06281)}
\end{equation}
which has one unstable pole in $p = 0.1081$ and one nonminimum phase zero in $z = 10$. We specify a maximum gain, phase and delay uncertainty of $\kappa = 1.6$, $\varphi = \pi/20$, and $\mathtt{t} = 1.5$, respectively. We use the method from Section~\ref{sec:approx_method}, with $T_0(s) \equiv 0$, to investigate if there is a stabilizing controller, and compute the interpolation values via numerical integration of \eqref{eq:outerrepr}.
Since the obtained Pick matrix is positive definite, there is a solution to \eqref{eq:T_tilde_shifted}. In fact, there are infinitely many rational solutions \cite{byrnes2001ageneralized}, \cite{fanizza2007passivity}, and among these we compute the so called maximum entropy solution \cite{georgio2003kullback}, \cite{blomqvist2005optimization}. The obtained complementary sensitivity function, as well as the open-loop transfer function and the controller, are given in \eqref{eq:controller_design_ex}.%
\footnote{Note that the controller \eqref{eq:K_controller_design_ex} has a pole in $-0.1659$, where the plant \eqref{eq:ex_plant} has a zero. Moreover, the controller has zeros in $-0.1490 \pm i0.2015$, where the plant has poles. That is why the degree of the controller is higher than the degree of the open-loop system.}
Finally, a Nyquist plot of the system, together with the forbidden region $\Gamma(i\omega)$,
is shown in Fig.~\ref{fig:num_ex_nyquist}.

\begin{figure*}[th]
\footnotesize
%
\hrulefill
\begin{subequations}\label{eq:controller_design_ex}
\begin{align}
& T(s) =
\frac{-7.576 s^{10} + 47.93 s^9 + 171.6 s^8 + 881.1 s^7 + 1574 s^6 + 2481 s^5 + 2322 s^4 + 1601 s^3 + 731.5 s^2 + 194.2 s + 29.03}
{s^{11} + 33.88 s^{10} + 321.5 s^9 + 1108 s^8 + 3230 s^7 + 4781 s^6 + 6208 s^5 + 4585 s^4 + 2787 s^3 + 957 s^2 + 223 s + 21.41}
\label{eq:T_controller_design_ex} \\
& P(s)K(s) =
\frac{-7.576 s^{10} + 47.93 s^9 + 171.6 s^8 + 881.1 s^7 + 1574 s^6 + 2481 s^5 + 2322 s^4 + 1601 s^3 + 731.5 s^2 + 194.2 s + 29.03}
{s^{11} + 41.46 s^{10} + 273.5 s^9 + 935.9 s^8 + 2348 s^7 + 3207 s^6 + 3727 s^5 + 2262 s^4 + 1186 s^3 + 225.5 s^2 + 28.78 s - 7.614}
\label{eq:PK_controller_design_ex}\\
& K(s) =
\frac{-75.76 s^{11} - 300.8 s^{10} - 1154 s^9 - 2182 s^8 - 3347 s^7 - 3420 s^6 - 2590 s^5 - 1407 s^4 - 526.3 s^3 - 135 s^2 - 21.03 s - 1.823}
{0.1s^{11} + 4.173 s^{10} + 28.49 s^9 + 101.2 s^8 + 261.3 s^7 + 387.9 s^6 + 467.9 s^5 + 338.7 s^4 + 192.8 s^3 + 63.07 s^2 + 13.44 s + 1.169}
\label{eq:K_controller_design_ex}
\end{align}
\end{subequations}
\hrulefill
\vspace*{4pt}
\end{figure*}

\section{Conclusions and future directions}\label{sec:conclusions}
In this work we show that robust stabilization can, in certain cases, be achieved by solving a non-classically constrained analytic interpolation problem. 
Moreover, by considering a sufficient condition for when this problem can be solved, we derive an algorithm that can be used to solve the interpolation problem and hence which can be used for controller synthesis.
In future work, it would be interesting to investigate if one can find necessary and sufficient conditions for solvability of  the non-classically constrained analytic interpolation problem.
Moreover, the interpolation problem is derived under assumptions on the set of disturbances, see Section~\ref{sec:ass}. In future work, we also want to investigate generalizations obtained by relaxing or changing some of the assumptions, for example by considering perturbations in the coprime factors of the plant. 


\appendix

\subsection{Proof of Theorem~\ref{thm:dist_to_set}}\label{app:proof_of_thm:thm:dist_to_set}
The theorem is proved using a sequence of lemmas. These in turn are proved under the assumption that the conditions in the theorem hold, i.e., that $\Omega$ satisfies the assumptions in Section~\ref{sec:ass}, that $\Lambda$ is given as in \eqref{eq:lambda}, that
 $T \in \HinfC \cap \mathcal{C}(i\mR)$, that $T(\infty)$ is well-defined, and that $T(\infty) \not \in \Lambda(\infty)$. Note that the latter is equivalent to $\dist(\Lambda(\infty), T(\infty)) > 0$.

\begin{lemma}\label{lem:encirclement}
For a fixed $\Delta \in \Omega$, $1 + T(s)(\Delta(s) - 1) \neq 0$ for all $s \in \mC_+$ if and only if the curve $\{ T(i \omega)(\Delta(i\omega) - 1) \mid \omega \in \mR \}$ does not intersect or encircle $-1$.
\end{lemma}

\hspace{-22pt}
\begin{proof}
This follows from the argument principle, see, e.g.,
\cite[Thm.~9.3]{astrom2008feedback}, and can be shown by using the same arguments as in the proof of the Nyquist stability criterion, see, e.g., \cite[Sec.~9.2]{astrom2008feedback}. In particular, observe that $f(s) := 1 + T(s)(\Delta(i\omega) - 1) \in \HinfC$, since $T, \Delta \in \HinfC$. This means that $f$ has no poles in $\tilde{\mC}_+$, and hence any encirclement of $-1$ will be clockwise and correspond to a zero of $f$ in the right half plane.
\end{proof}

\begin{lemma}\label{lem:not_zero_s_omega}
The following two statements are equivalent:
\begin{enumerate}
\item[i)] $1 + T(s)(\Delta(s) - 1) \neq 0$ for all $s \in \tilde{\mC}_+$ and $\Delta \in \Omega$,
\item[ii)] $1 + T(i\omega)(\Delta(i\omega) - 1) \neq 0$ for all $\omega \in \mR$ and $\Delta \in \Omega$.
\end{enumerate}
\end{lemma}

\hspace{-22pt}
\begin{proof}
We proceed by proving the equivalence between the negations of the two statements, i.e., the equivalence between the two statements
\begin{enumerate}
\item[i\ensuremath{'})] there exists a $\Delta \in \Omega$ and an $s \in \tilde{\mC}_+$ such that \newline $1 + T(s)(\Delta(s) - 1) = 0$,
\item[ii\ensuremath{'})] there exists a $\Delta \in \Omega$ and an $\omega \in \mR$ such that \newline $1 + T(i\omega)(\Delta(i\omega) - 1) = 0$.
\end{enumerate}
That the second statement implies the first is trivial, since $i\mR \subset \tilde{\mC}_+$, and we can thus take the same point and the same $\Delta$. Moreover, if in the first statement the point $s$ is such that $s \in i\mR$, then in same way the first statement implies the second. The lemma thus follows by showing that, if $i\ensuremath{'})$ holds for $s \in \mC_+$, then it implies $ii\ensuremath{'})$, and that neither $i\ensuremath{'})$ nor $ii\ensuremath{'})$ can hold for $s = \infty$, i.e., that there can be no sequence $(s_n)_n \subset \bar{\mC}_+$ such that $|s_n| \to \infty$ as $n \to \infty$ and such that
\begin{equation}\label{eq:limsup_proof}
\lim_{n \to \infty} 1 + T(s_n)(\Delta(s_n) - 1) = 0.
\end{equation}

To show the latter, note that in points where $\Delta(s_n) = 1$ we have that $1 + T(s_n)(\Delta(s_n) - 1) = 1$ and hence such points can be removed from any potential sequence where $0$ is the limit. Without loss of generality, we therefore only consider sequences $(s_n)_n$ such that $ \Delta(s_n) \neq 1$ for all $n$. For such sequences, \eqref{eq:limsup_proof} is true if and only if 
\[
\lim_{n \to \infty} \left| \frac{1}{\Delta(s_n) - 1} + T(s_n) \right| = 0.
\]
However, this is not possible since by assumption $T(\infty) \not \in \Lambda(\infty)$.

In the case $i\ensuremath{'})$ holds for some $\Delta \in \Omega$ and $s \in \mC_+$, by Lemma~\ref{lem:encirclement} we have that the curve $\{ T(i \omega)(\Delta(i \omega) - 1) \mid \omega \in \mR \}$ intersects or encircles $-1$ at least once. However, the curve varies continuously pointwise with the disturbance and for the disturbance $I \in \Omega$ the curve is identically $0$. Therefore, there must exist a $\tilde{\Delta}$ and an $\tilde{\omega} \in \mR$ such that $T(i \tilde{\omega})(\tilde{\Delta}(i\tilde{\omega}) - 1) = -1$, which shows that i') implies ii') in this case. This completes the proof.
\end{proof}

To prove the theorem, we first prove that $1)$ implies $2)$. To this end, assume that there is an $\varepsilon > 0$ such that 
\begin{equation}\label{eq:proof_eq1}
|1 + T(s)(\Delta(s) - 1)| = |1 - T(s)(1 - \Delta(s))| \geq \varepsilon
\end{equation}
for all $s \in \tilde{\mC}_+$ and all $\Delta \in \Omega$.
Now, fix $\Delta$ and assume that $\Delta(\tilde{s}) = 1$. By continuity of $\Delta$ we can construct a region $D_{\tilde{s}}^{\Delta}$ around $\tilde{s}$ such that $|1/(1 - \Delta(s))| \geq 2M$ for all $s \in D_{\tilde{s}}$, where $M := \| T \|_{\Hinf}$. Therefore, by the reverse triangle inequality, for all $s \in D_{\tilde{s}}^{\Delta}$ we have
\[
\left| \frac{1}{1 - \Delta(s)} - T(s) \right| \geq \left| \left| \frac{1}{1 - \Delta(s)} \right| -  |T(s)| \right| \geq M.
\]
Let $D^{\Delta} = \cup_{s \text{ such that } \Delta(s) = 1} D_s^{\Delta} $. On $\tilde{\mC}_+ \setminus D^{\Delta}$ we have $\Delta(s)$ bounded away from $1$ and thus \eqref{eq:proof_eq1} implies that 
\[
\left| \frac{1}{1 - \Delta(s)} - T(s) \right| \geq \varepsilon \left| \frac{1}{1 - \Delta(s)} \right|.
\]
Since
$\sup_{\Delta \in \Omega} \| \Delta \|_{\Hinf} =: N < \infty$, and since $\inf_x 1/f(x) = 1/\sup_x f(x)$ for nonnegative functions $f$, we thus have that
\begin{align*}
& \inf_{\substack{s \in \bar{\mC}_+ \\ \Delta \in \Omega}} \! \left| \frac{1}{1 \! - \! \Delta(s)} \right|
= \frac{1}{\sup_{\!\substack{s \in \bar{\mC}_+ \\ \Delta \in \Omega}} \! \left| 1 \! - \! \Delta(s) \right|}
\geq \frac{1}{1 \! + \! \sup_{\!\substack{s \in \bar{\mC}_+ \\ \Delta \in \Omega}} \! | \Delta(s) |}  \\
& \phantom{xxx} = \frac{1}{1 + \sup_{\Delta \in \Omega} \|\Delta(s)\|_{\Hinf}} 
= \frac{1}{1 + N}.
\end{align*}
For a fixed $\Delta$ we therefore have
\[
\inf_{s \in \bar{\mC}_+ \! \setminus D^{\Delta}} \left| \frac{1}{1 \! - \! \Delta(s)} - T(s)  \right| \! \geq \!\!\! \inf_{s \in \bar{\mC}_+ \! \setminus D^{\Delta}} \varepsilon \! \left| \frac{1}{1 \! - \! \Delta(s)} \right| \! \geq \! \frac{\varepsilon}{1 \! + \! N},
\]
and hence
\begin{align*}
& \inf_{s \in \bar{\mC}_+} \left| \frac{1}{1 - \Delta(s)} - T(s) \right| = \\
&  \min
\left\{
\! \inf_{s \in D^{\Delta}} \! \left| \frac{1}{1 \! - \! \Delta(s)} - T(s) \right|
\! , \!
\inf_{s \in \bar{\mC}_+ \! \setminus D^{\Delta}} \! \left| \frac{1}{1 \! - \! \Delta(s)} - T(s) \right|
\right\} \\
& \geq \min \{ M, \epsilon /(1 + N) \} =: \tilde{\varepsilon} > 0. 
\end{align*}
Finally, this gives the sought inequality
\begin{align*}
& \inf_{\omega \in \mR} \dist(\Lambda(i\omega), T(i\omega)) \! = \! \inf_{\omega \in \mR} \inf_{\Delta \in \Omega} \left| \frac{1}{1 - \Delta(i\omega)} - T(i\omega) \right| \\
& \phantom{xxx} \geq \inf_{s \in \bar{\mC}_+} \inf_{\Delta \in \Omega} \left| \frac{1}{1 - \Delta(s)} - T(s) \right| \\
& \phantom{xxx} = \inf_{\Delta \in \Omega} \inf_{s \in \bar{\mC}_+} \left| \frac{1}{1 - \Delta(s)} - T(s) \right| \\
& \phantom{xxx} \geq \inf_{\Delta \in \Omega} \min \{ M, \epsilon/(1 + N) \}  = \inf_{\Delta \in \Omega} \tilde{\varepsilon} = \tilde{\varepsilon} > 0.
\end{align*}

To prove the converse, assume there is an $\varepsilon > 0$ such that
\[
\inf_{\omega \in \mR} \dist(\Lambda(i\omega), T(i\omega)) \geq \varepsilon.
\]
In particular, this means that for all $\omega \in \mR$ and all $\Delta \in \Omega$, $| 1/(1 - \Delta(i\omega)) - T(i\omega) | \neq 0$, and thus that
\begin{equation}\label{eq:proof_thm_dist_cut}
| 1 - T(i\omega)(1 - \Delta(i\omega)) | \neq 0.
\end{equation}
By Lemma~\ref{lem:not_zero_s_omega}, this means that $1 + T(s)(\Delta(s) - 1) \neq 0$ for all $s \in \tilde{\mC}_+$ and $\Delta \in \Omega$. Left to show is that this implies that
\[
\inf_{\substack{s \in \bar{\mC}_+ \\ \Delta \in \Omega}} |1 + T(s)(\Delta(s) - 1)| > 0.
\] 
To this end, assume that the latter is false. Then there is a sequence $(s_n, \Delta_n)_n \subset \bar{\mC}_+ \times \Omega$ such that
\begin{equation}\label{eq:contradiction_statement_proof}
\lim_{n \to \infty} 1 + T(s_n)(\Delta_n(s_n) - 1) = 0.
\end{equation}
First assume that the sequence $(s_n)_n$ is bounded. Then there is a subsequence converging to a point $s_\infty \in \bar{\mC}_+$.  Since $\Omega$ is bounded,
by assumption there is a sub-subsequence of $(\Delta_n)_n$ that converges pointwise to some $\Delta_\infty \in \Omega$, and thus for this sub-subsequence we have that $\lim_{n \to \infty} \Delta_n(s_n) = \Delta_\infty (s_\infty)$. However, by continuity this means that $1 + T(s_\infty)(\Delta_\infty (s_\infty) - 1) = 0$, which contradicts \eqref{eq:proof_thm_dist_cut} and Lemma~\ref{lem:not_zero_s_omega}. Thus any such sequence $(s_n)_n$ must be unbounded.
Next, note that we must have $\lim_{n \to \infty} \Delta_n(s_n) \neq 1$ since $T(\infty)$ is assumed to be well-defined and bounded. Therefore, we can equivalently consider the limit
\[
\lim_{n \to \infty} \left| \frac{1}{\Delta_n(s_n) - 1} + T(s_n) \right| = 0.
\]
However, by definition that would mean that $T(\infty) \in \Lambda(\infty)$, which by assumption is false. Therefore, the infimum is strictly larger than 0, which proves the statement.
\qed

\subsection{Proof of Lemma~\ref{lem:regularity_phi}}\label{app:proof_lem:regularity_phi}
Clearly, $\phi$ is nonnegative. Moreover,
\begin{align*}
\phi(i\omega) &= \sup_{\Delta \in \Omega} |\Delta(i\omega) - 1| \leq \sup_{\Delta \in \Omega} |\Delta(i\omega)| + 1 \\ 
&\leq \sup_{\Delta \in \Omega} \|\Delta \|_{\HinfC} + 1 = N +1.
\end{align*}
Finally, for each $\Delta \in \Omega$ the function $\omega \mapsto |\Delta(i\omega) - 1|$ is a continuous function from $\mR$ to $\mR$, and thus lower semicontinuous. The result thus follows since the supremum over any collection of lower semicontinuous functions is lower semicontinuous (see, e.g., \cite[p.~38]{rudin1987real}).
\qed


\balance

\bibliographystyle{plain}
\bibliography{ref}

\begin{thebibliography}{10}

\bibitem{blomqvist2005optimization}
A.~Blomqvist and R.~Nagamune.
\newblock Optimization-based computation of analytic interpolants of bounded
  complexity.
\newblock {\em Systems \& Control Letters}, 54(9):855--864, 2005.

\bibitem{byrnes2001ageneralized}
C.I. Byrnes, T.T. Georgiou, and A.~Lindquist.
\newblock A generalized entropy criterion for {Nevanlinna-Pick} interpolation
  with degree constraint.
\newblock {\em IEEE Transactions on Automatic Control}, 46(6):822--839, 2001.

\bibitem{callier1978algebra}
F.M. Callier and C.A. Desoer.
\newblock An algebra of transfer functions for distributed linear
  time-invariant systems.
\newblock {\em IEEE Transactions on Circuits and Systems}, 25(9):651--662,
  1978.

\bibitem{chen1982necessary}
M.-J. Chen and C.A. Desoer.
\newblock Necessary and sufficient condition for robust stability of linear
  distributed feedback systems.
\newblock {\em International Journal of Control}, 35(2):255--267, 1982.

\bibitem{desoer1980generalized}
C.A. Desoer and Y.-T. Wang.
\newblock On the generalized {N}yquist stability criterion.
\newblock {\em IEEE Transactions on Automatic Control}, 25(2):187--196, 1980.

\bibitem{doyle1981multivariable}
J.~Doyle and G.~Stein.
\newblock Multivariable feedback design: Concepts for a classical/modern
  synthesis.
\newblock {\em IEEE transactions on Automatic Control}, 26(1):4--16, 1981.

\bibitem{doyle1992feedback}
J.C. Doyle, B.A. Francis, and A.~Tannenbaum.
\newblock {\em Feedback control theory}.
\newblock Macmillan, New York, N.Y., 1992.

\bibitem{fanizza2007passivity}
G.~Fanizza, J.~Karlsson, A.~Lindquist, and R.~Nagamune.
\newblock Passivity-preserving model reduction by analytic interpolation.
\newblock {\em Linear Algebra and its Applications}, 425(2-3):608--633, 2007.

\bibitem{foias1996robust}
C.~Foias, H.~{\"O}zbay, and A.~Tannenbaum.
\newblock {\em Robust control of infinite dimensional systems}.
\newblock Springer, Berlin, Heidelberg, 1996.

\bibitem{georgio2003kullback}
T.T. Georgiou and A.~Lindquist.
\newblock {K}ullback-{L}eibler approximation of spectral density functions.
\newblock {\em IEEE Transactions on Information Theory}, 49(11):2910--2917,
  2003.

\bibitem{helton1998classical}
J.W. Helton and O.~Merino.
\newblock {\em Classical Control Using {H}$^\infty$ Methods: Theory,
  Optimization, and Design}.
\newblock SIAM, Philadelphia, P.A., 1998.

\bibitem{hoffman1962banach}
K.~Hoffman.
\newblock {\em {B}anach Spaces of Analytic Functions}.
\newblock Prentice-Hall, Englewood Cliffs, N.J., 1962.

\bibitem{ju2016further}
P.~Ju and H.~Zhang.
\newblock Further results on the achievable delay margin using {LTI} control.
\newblock {\em IEEE Transactions on Automatic Control}, 61(10):3134--3139,
  2016.

\bibitem{khargonekar1986robust}
P.P. Khargonekar and K.~Poolla.
\newblock Robust stabilization of distributed systems.
\newblock {\em Automatica}, 22(1):77--84, 1986.

\bibitem{khargonekar1985non}
P.P. Khargonekar and A.~Tannenbaum.
\newblock Non-{E}uclidian metrics and the robust stabilization of systems with
  parameter uncertainty.
\newblock {\em IEEE Transactions on Automatic Control}, 30(10):1005--1013,
  1985.

\bibitem{kimura1984robust}
H.~Kimura.
\newblock Robust stabilizability for a class of transfer functions.
\newblock {\em IEEE Transactions on Automatic Control}, 29(9):788--793, 1984.

\bibitem{middleton2007achievable}
R.H. Middleton and D.E. Miller.
\newblock On the achievable delay margin using {LTI} control for unstable
  plants.
\newblock {\em IEEE Transactions on Automatic Control}, 52(7):1194--1207, 2007.

\bibitem{qi2014fundamental}
T.~Qi, J.~Zhu, and J.~Chen.
\newblock Fundamental bounds on delay margin: When is a delay system
  stabilizable?
\newblock In {\em 33rd Chinese Control Conference (CCC)}, pages 6006--6013.
  IEEE, 2014.

\bibitem{qi2017fundamental}
T.~Qi, J.~Zhu, and J.~Chen.
\newblock Fundamental limits on uncertain delays: {W}hen is a delay system
  stabilizable by {LTI} controllers?
\newblock {\em IEEE Transactions on Automatic Control}, 62(3):1314--1328, 2017.

\bibitem{astrom2008feedback}
K.J. \r{A}str\"{o}m and R.M. Murray.
\newblock {\em Feedback systems}.
\newblock Princeton university press, Princeton, N.J., 2008.

\bibitem{ringh2018lower}
A.~Ringh, J.~Karlsson, and A.~Lindquist.
\newblock Lower bounds on the maximum delay margin by analytic interpolation.
\newblock In {\em 2018 IEEE Conference on Decision and Control (CDC)}, pages
  5463--5469. IEEE, 2018.

\bibitem{ringh2019ananalytic}
A.~Ringh, J.~Karlsson, and A.~Lindquist.
\newblock An analytic interpolation approach to stability margins with emphasis
  on time delay.
\newblock {\em arXiv:1912.08734}, 2019.

\bibitem{rudin1987real}
W.~Rudin.
\newblock {\em Real and {C}omplex {A}nalysis}.
\newblock McGraw-Hill, New York, 3d edition, 1987.

\bibitem{schiff1993normal}
J.L. Schiff.
\newblock {\em Normal Families}.
\newblock Springer, New York, N.Y., 1993.

\bibitem{tannenbaum1980feedback}
A.~Tannenbaum.
\newblock Feedback stabilization of linear dynamical plants with uncertainty in
  the gain factor.
\newblock {\em International Journal of Control}, 32(1):1--16, 1980.

\bibitem{youla1974single}
D.C. Youla, J.J. Bongiorno~Jr, and C.N. Lu.
\newblock Single-loop feedback-stabilization of linear multivariable dynamical
  plants.
\newblock {\em Automatica}, 10(2):159--173, 1974.

\bibitem{zames1983feedback}
G.~Zames and B.~Francis.
\newblock Feedback, minimax sensitivity, and optimal robustness.
\newblock {\em IEEE Transactions on Automatic Control}, 28(5):585--601, 1983.

\bibitem{zhou1996robust}
K.~Zhou, J.C. Doyle, and K.~Glover.
\newblock {\em Robust and optimal control}.
\newblock Prentice-Hall, Englewood Cliffs, N.J., 1996.

\end{thebibliography}

\end{document}